\documentclass{ws-ijbc}
\usepackage{ws-rotating}     
\usepackage{graphicx}
\usepackage{epstopdf}
\begin{document}

\catchline{}{}{}{}{} 

\markboth{Zhang}{weighted m-clique annex operation}

\title{Spectral Analysis and its applications for a class of scale-free network based on the weighted m-clique annex operation}

\author{Zhizhuo Zhang}
\address{School of Mathematics, Southeast University,\\ Nanjing 211189, PR China\\
zhizhuo\_zhang@163.com}

\author{Jinde Cao\footnote{Corresponding Author.}}
\address{School of Mathematics, Southeast University, \\
Nanjing 211189, PR China\\
Yonsei Frontier Lab, Yonsei University,\\
Seoul 03722, South Korea\\
jdcao@seu.edu.cn}

\author{Bo Wu}%
\address{ School of Applied Mathematics, Nanjing University of Finance and Economics,\\
 Nanjing, 210023, P.R. China\\
bowu8800@nufe.edu.cn}%

\author{Ardak Kashkynbayev}

\address{Department of Mathematics, Nazarbayev University,\\
Nur-Sultan 010000, Kazakhstan\\
ardak.kashkynbayev@nu.edu.kz}

\maketitle

\begin{abstract}
The spectrum of network is an important tool to study the function and dynamic properties of network, and graph operation and product is an effective mechanism to construct a specific local and global topological structure. In this study, a class of weighted $m-$clique annex operation $\tau_m^r(\cdot)$ controlled by scale factor $m$ and weight factor $r$ is defined, through which an iterative weighted network model $G_t$ with small-world and scale-free properties is constructed. In particular, when the number of iterations $t$ tends to infinity, the network has transfinite fractal property. Then, through the iterative features of the network structure, the iterative relationship of the eigenvalues of the normalized Laplacian matrix corresponding to the network is studied. Accordingly, some applications of the spectrum of the network, including the Kenemy constant, Multiplicative Degree-Kirchhoff index and the number of weighted spanning trees, are further given. In addition, we also study the effect of the two factors controlling network operation on the structure and function of the iterative weighted network $G_t$, so that the network operation can better simulate the real network and have more application potential in the field of artificial network.
\end{abstract}


\section{Introduction}
As a new and interdisciplinary science, complex network has attracted the attention of many scholars because of its application in various fields\cite{newman2003structure}. Modern science has confirmed that theory of graph spectra is an important tool for studying complex networks\cite{wu2019spectral,yi2018scale,qi2018extended,zeng2021spectra,van2010graph}. Many functions or dynamic properties on a complex network are difficult to directly respond to the topology of the network, but can be studied through the properties of the spectrum, such as: mixing time, Laplacian energy, hitting time, Kemeny constant, Multiplicative Degree-Kirchhoff index, the number of spanning trees, etc\cite{mehatari2015effect,julaiti2013eigenvalues,tetali1991random,chen2007resistance,chang2014spanning,qi2015novel}.  Moreover, studies of complex networks in real world show that these network models often have some common overall properties, such as small-world properties, scale-free properties, and so on\cite{watts1998collective,barabasi1999emergence}. Therefore, it has become a hot topic in network science to study the dynamic properties of complex networks with this properties by using the spectrum of networks.

However, for large-scale random complex networks, it is often very difficult to study their spectra due to the complexity of their corresponding matrices. Therefore, considerable attention has been used to find the generation mechanism and model of the network showing the salient features of the real network, and to study the properties of the spectrum on it\cite{pralat2011edge,zhang2011farey,barrat2004weighted,barrat2004modeling}. Graph operation and product is a representative of these network generation mechanisms, which has been widely used to generate regular network models with scale-free, small-world, and self-similarity properties. Furthermore, given that complex networks often have some special local structures, such as communities, motifs and cliques\cite{girvan2002community,milo2002network,tsourakakis2015k}, and graph operations and products can iteratively generate a large scale network with some special local structures from a small initial network, more and more attention has been paid to this generating mechanism of the complex network. In addition, due to the certainty of network operation rules, the properties of the spectrum are often studied systematically, which cannot be achieved in the real complex network\cite{wu2019spectral,yi2018scale,qi2018extended,zeng2021spectra}. At present, many network operations have been used to design and construct complex network models and to study the topological and dynamic properties on them, including: $q-$subdivision, planar triangulation, Cartesian product, hierarchical product, corona product, Kronecker product, etc.\cite{wood2005acyclic,jin2017maximum,imrich2000product,barriere2016deterministic,sharma2017structural,mahdian2007stochastic}. However, the network iteration methods in the above articles all use nodes as the basic unit, and the iterative model with local node sets as the unit has never been studied.

In this paper, we define a class of weighted $m-$clique annex operation, in which we attach an $m-$complete subgraph to each edge of the network. Therefore, the network model constructed by this operation has significant local characteristics. In addition, the weighted network model constructed by iterative operation has significant small world, scale-free and self-similarity. Moreover, when the number of iterations $t$ tends to infinity, the network has transfinite fractal property. Then, based on the change of network topology structure by the weighted $m-$clique annex operation, the iterative relation of spectrum (eigenvalue) of normalized Laplacian matrix will be analysed. Furthermore, we will calculate the analytic expressions of three characteristic quantities on the weighted network model, which are Kemeny constant, Multiplicative Degree-Kirchhoff index and the number of weighted spanning trees. These characteristics are closely related to the overall topological or dynamic properties of the network, which can be used not only to analyze the influence of the special topological structure of the network on the network function and its topological properties, but also have application potential in the field of science and engineering. It is worth mentioning that the weighted $m-$clique annex operation is controlled by the scale factor $m$ and the weight factor $r$, in which the change of the scale factor $m$ will greatly affect the overall topology structure of the network, while the weight factor $r$ can affect the network functions and topological properties without changing the network structure. Therefore, the combination of the two factors can not only better simulate the real network model, but also provide better controllability for some artificial network systems.

This article will be divided into the following four parts: in the second Section, we will introduce the specific definition of the weighted $m-$clique annex operation and some topological properties of the iterative weighted network constructed from it, including the number of nodes, degree distribution, and transfractal dimension. In the third Section, the iterative relationship of the normalized Laplaacian matrix from the changes in the network topology caused by network operations will be derived, and the iterative expression of its spectrum will be proved. In the fourth Section, three applications of network normalized Laplacian matrix spectrum will be presented. In addition, the relationship between Kenemy constants and scale factors and weight factors will be numerically simulated to analyze the impact of the above two factors on network functions and topological properties. Finally, in fifth Section, we will summarize the whole article.

\section{Properties of weighted $m-$clique annex operation}

\subsection{Construction method of the operation}
In this Section, the construction of weighted $m-$clique annex operation and some topological properties are introduced . First of all, it should be clear that the $m-$clique mentioned here is a complete subgraph with $m$ nodes, denoted as $K_m$. Since the initial network is only required to be a connected network, which is denoted as $G_0$, the weighted $m-$clique annex operation has strong applicability. Then, the construction of the weighted $m-$clique annex operation is as follows: for any edge with weight of $\omega$ in network $G_0$, add a clique $K_m$, and connect all nodes in the $K_m$ to the two nodes connected by the edge, and then specify that the weight of all newly generated edges is $r\omega$. The new weighted network generated by the network operation of the initial network $G_0$ is denoted as $G_1$. The weighted $m-$clique annex operation is controlled by and only by two parameters, namely scale factor $m(m\in\mathbb{N}^+)$ and weight factor $r(r\in\mathbb{R}^+)$. Therefore, the network operation process is denoted as $\tau^r_m(\cdot)$, so the following relation can be obtained:
$$
G_1=\tau^r_m(G_0).
$$
Obviously, in this operation, each edge in $G_0$ will correspond to a weighted clique $K_m$, and the weight of the edge connected by all nodes in this clique is determined by the edge in this $G_0$, so we call this edge the "parent" of this weighted clique $K_m$. For example, when the size factor $m=4$, the process of network operation $\tau^r_4(\cdot)$ is shown in Fig.1.
\begin{figure}[t]
\centering
\includegraphics[width=12cm]{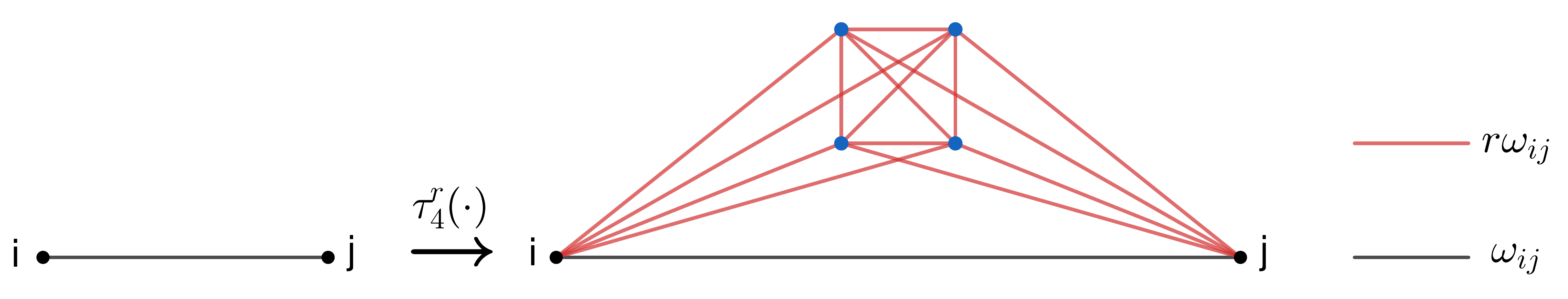}
\caption{Schematic diagram of the weighted $m-$clique annex operation $\tau^r_m(\cdot)$ with $m=4$: black nodes $i$ and $j$ are the nodes existing in the original network, and the weight of the edge connecting the two nodes is denoted as $\omega_{ij}$; the blue nodes are the new nodes generated after network operation $\tau^r_4(\cdot)$, where the weight of the red edge is $r\omega_{ij}$.}
\end{figure}

Naturally, for any initial network $G_0$, the weighted $m-$clique annex operation $\tau^r_m(\cdot)$ can be repeatedly used to generate a complex network model, denoted as $G_t$, where $t$ represents the number of iterations of the operation. Therefore, network $G_t$ can be expressed as follows:
$$
G_t=\underbrace{\tau^r_m(\tau^r_m(\cdots\tau^r_m(G_0)\cdots))}_{t}.
$$
According to the definition of operation $\tau^r_m(\cdot)$, as the number of iterations increases, the size of the network will increase exponentially, but the diameter of the network will only increase linearly, so the iterative weighted network has a small-world property.
In order to facilitate later calculation, when we consider the relevant properties of weighted iterative network $G_t$, each edge in the initial network $G_0$ is set as unit weight.

\begin{remark}
In the network operator $\tau^r_m(\cdot)$, the introduction of the local weighted clique $K_m$ has both practical and theoretical significance. Firstly, the reality is based on the fact that during the dynamic growth of the actual network, the newly generated nodes are to be connected to the original network not necessarily in the form of isolated nodes, but based on a certain local structure as the basic unit. Therefore, the designed network operator $\tau^r_m(\cdot)$ can clearly reflect this feature. Second, by studying the dynamic characteristics of the network model generated by the network operator $\tau^r_m(\cdot)$, the extent to which the dynamic growth network $G_t$ is affected by such a local structure $K_m$ will be explained theoretically.
\end{remark}

\subsection{Topological properties of operation}

Then, the iterative rules of nodes, edges and strengths in the network under this operation will be studied. First, we specify that the symbol $|\cdot|$ represents the number of elements in the set. For any initial network $G_0$, the set of nodes and the set of edges are denoted as $V_0$ and $E_0$, respectively. Obviously, after the weighted $m-$clique annex operation $\tau^r_m(\cdot)$, the number of nodes and the number of edges in network $G_1$ satisfy the following relation:
$$
|E_1|=|E_0|\cdot(\frac12m^2+\frac32m+1)~~~\textrm{and}~~~ |V_1|=|V_0|+m\cdot|E_0|.
$$
Then, the sets of nodes and edges in the iterated weighted network $G_t$ are denoted as $N_t$ and $E_t$ respectively. From the above relationship, it can be naturally deduced that:
\begin{eqnarray}
|E_t|&=&|E_0|\cdot(\frac12m^2+\frac32m+1)^t\label{eq:2.2.1}\\
|V_t|&=&|V_0|+\frac{2|E_0|}{m+3}\Big[(\frac12m^2+\frac32m+1)^t-1\Big].\label{eq:2.2.2}
\end{eqnarray}

Then, the set of nodes already existing in the previous generation network $G_{t-1}$ contained in $V_t$ is denoted as $\hat V_t$. Therefore, the rest nodes are newly generated after the $t$-th operation $\tau_m^r(\cdot)$, and the set of this nodes is denoted as $\bar V_t$. Moreover, for any node $i\in V_t$, the set formed by nodes adjacent to node $i$ is denoted as $V^i_t$. Then, for any node $j\in V^i_t$, there must be an edge between node $i$ and node $j$, denoted as $e_{ij}$, and the weight of this edge denoted as $\omega_{ij}$. In addition, for the convenience of later calculation, the following provisions are made:
$$
\bar V_t^i=V_t^i\cap\bar V_t,~~ \hat V_t^i=V_t^i\cap\hat V_t,~~ \bar V_t^{ij}=V_t^i\cap V_t^j\cap\bar V_t.
$$
Then, for any node $i\in V_t$, its strength is defined as:
$$
s_t^i=\sum_{j\in V^i_t}\omega_{ij}.
$$
Naturally, if $i\in\hat V_t$, from the construction mode of operation $\tau_m^r(\cdot)$, it can be obtained that $s_t^i$ satisfies the following relation:
\begin{eqnarray}
s_t^i=(mr+1)s_{t-1}^i ~~~t\geq1.\label{eq:2.2.3}
\end{eqnarray}
If node $i$ belongs to set $\bar V_t$, the parent edge of this node is denoted as $e^*_i$, and the weight of this edge is denoted as $\omega^*_i$, then the strength of this node must satisfy:
\begin{eqnarray}
s_t^i=(m+1)r\omega^*_i.\label{eq:2.2.4}
\end{eqnarray}
Therefore, the total strength of the iterative weighted network $G_t$, denoted as $S_t$, can be defined as:
$$
S_t=\sum_{e_{ij}\in E_t}\omega_{ij}=\frac12\sum_{i\in V_t}s_t^i.
$$
The total strength $S_t$ of weighted iterative network $G_t$ can be obtained by network operation $\tau_m^r(\cdot)$, which satisfies the following relation:
$$
S_t=(\frac12rm^2+\frac32rm+1)S_{t-1}=(\frac12rm^2+\frac32rm+1)^tS_{0}.
$$
In the iterative weighted network $G_t$, since the weight of each edge of the initial network $G_0$ is defined as $1$, it can be obtained that $S_0=E_0$.

\subsection{Degree distribution of the iterative network $G_t$}

Next, the influence of the weighted $m-$clique annex operation $\tau^r_m(\cdot)$ on the network degree distribution will be discussed. Firstly, for any node $i\in V_t$, its degree, denoted as $d_t^i$, is defined as the number of adjacent nodes of the node $i$, which satisfy $d_t^i=|V_t^i|$. Similar to the node strength in network $G_t$, it is easy to prove that for any node $i\in V_t$, its degree $d_t^i$ must satisfy the following relation:
$$
d_t^i=\left\{
\begin{array}{lcl}
(m+1)d_{t-1}^i ~~& , & i\in \hat V_t \\
m+1 ~~& , & i\in \bar V_t
\end{array}
\right..
$$
From the above relation, it can be further obtained:
$$
d_t^i=\left\{
\begin{array}{lcl}
d_0^i(m+1)^t ~~& , & \textrm{if}~~ i\in V_0 \\
(m+1)^k ~~\textrm{and}~~ 1\leq k \leq t ~~& , & \textrm{other}
\end{array}
\right.,
$$
where $d_0^i$ represents the degree of node $i$ in network $G_0$.

Then, for the iterative network $G_t$, its degree distribution, denoted as $P^t(\alpha)$, is defined as the probability that a node is randomly selected from set $V_t$ with degree $\alpha$. However, since the degree distribution of the network is generally a discrete function, we pay more attention to the cumulative degree distribution on the iterative network $G_t$, denoted as $P^t_{cum}(\alpha)$, which is defined as the probability that the degree of randomly selected node is greater than or equal to $\alpha$, that is:
$$
P^t_{cum}(\alpha)=\sum_{x=\alpha}^\infty P^t(x).
$$

Since the network structure of the initial network $G_0$ is known, the degree distribution $P^0(\alpha)$ and the accumulative degree distribution $P^0_{cum}(\alpha)$ on the network are known initial conditions. Therefore, for node $i\in V_0$ in iterative network $G_t$, its degree distribution must satisfy the following conditions:
\begin{eqnarray}
P^t(d_t^i)=P^t(d_0^i(m+1)^t)=\frac{|V_0|}{|V_t|}P^0(d_0^i)=P^0(d_0^i)\cdot \Big(1+\frac{2|E_0|}{|V_0|(m+3)}\Big[(\frac12m^2+\frac32m+1)^t-1\Big]\Big)^{-1}\label{eq:2.3.1}
\end{eqnarray}
According to the weighted $m-$clique annex operation $\tau^r_m(\cdot)$, the degree of nodes generated in the same iteration operation in the iterative network $G_t$ is equal, and the degree of nodes generated in different generations must not be equal. Therefore, based on the iterative rule of node degree, it is easy to prove that the number of nodes with degree $(m+1)^k$ belonging to set $V_t/V_0$ is $E_0\cdot m(\frac12m^2+\frac32m+1)^{t-k}$, where $1\leq k\leq t$. Therefore, for node $i\in V_t/V_0$ with degree $(m+1)^k$ in network $G_t$, its degree distribution satisfies the following relation:
\begin{eqnarray}
P^t(d_t^i)\!=\!P^t((m+1)^k)\!=\!\frac{E_0\cdot m(\frac12m^2+\frac32m+1)^{t-k}}{|V_0|+\frac{2|E_0|}{m+3}\big[(\frac12m^2+\frac32m+1)^t-1\big]}.\label{eq:2.3.2}
\end{eqnarray}
From Eq.(\ref{eq:2.3.1}) and Eq.(\ref{eq:2.3.2}), the degree distribution on the iterative network $G_t$ can be obtained, and then its cumulative degree distribution can be solved.

If $\alpha\in\mathbb{R}^+$ and $\alpha\leq(m+1)^t$, according to the definition, the cumulative degree distribution $P^t_{cum}(\alpha)$ on the iterative network $G_t$ can be expressed as:
\begin{eqnarray}\label{eq:2.3.3}
P^t_{cum}(\alpha)&=&1-\sum_{k=1}^\beta P^t((m+1)^k)\nonumber\\
&=&1-\sum_{k=1}^{\beta}\frac{E_0\cdot m(\frac12m^2+\frac32m+1)^{t-k}}{|V_0|+\frac{2|E_0|}{m+3}\big[(\frac12m^2+\frac32m+1)^t-1\big]}\nonumber\\
&=&\frac{|V_0|(m+3)+2|E_0|[(\frac12m^2+\frac32m+1)^{t-\beta}-1]}{|V_0|(m+3)+2|E_0|[(\frac12m^2+\frac32m+1)^t-1]},
\end{eqnarray}
where $\beta$ is the largest positive integer that satisfy $(m+1)^\beta<\alpha$, so it can be obtained that $\beta=\lceil\frac{\ln\alpha}{\ln(m+1)}\rceil-1$.

If $\alpha\in\mathbb{R}^+$ and $\alpha>(m+1)^t$, it is easy to deduce from Eq.(\ref{eq:2.3.1}) that the accumulative degree distribution $P_{cum}^t(\alpha)$ on the iterative network $G_t$ satisfies the following relation:
\begin{eqnarray}\label{eq:2.3.4}
&&P^t_{cum}(\alpha)=\frac{|V_0|}{|V_t|}P^0_{cum}(\beta)\nonumber=P^0_{cum}(\beta) \Big(1+\frac{2|E_0|}{|V_0|(m+3)}\Big[(\frac12m^2+\frac32m+1)^t-1\Big]\Big)^{-1},
\end{eqnarray}
where $\beta=\lceil\frac{\alpha}{(m+1)^t}\rceil$. Therefore, from Eq.(\ref{eq:2.3.3}) and Eq.(\ref{eq:2.3.4}), the analytic expression of the accumulative degree distribution of iterated network $G_t$ have obtained.

When the number of iterations in the network $G_t$ is large enough so that $|V_t|\gg |V_0|$, we can obtain:
$$
P^t_{cum}(\alpha)=\frac{|V_0|}{|V_t|}P^0_{cum}(\beta)\approx0 ,~~ \alpha>(m+1)^t.
$$
Therefore, the accumulative degree distribution $P^t_{cum}(\alpha)$ on network $G_t$ can be approximately expressed as:
$$
P^t_{cum}(\alpha)\approx\frac{|V_0|(m+3)+2|E_0|[\alpha^{1-\gamma}(\frac12m^2+\frac32m+1)^t-1]}
{|V_0|(m+3)+2|E_0|[(\frac12m^2+\frac32m+1)^t-1]},
$$
where $\gamma=2+\frac{\ln(\frac{m}{2}+1)}{\ln(m+1)}$. Naturally, for large number of iterations $t$, it can be proofed that:
$$
P^t_{cum}(\alpha) \sim \alpha^{1-\gamma} ~~~\textrm{and}~~~ P^t(\alpha) \sim \alpha^{-\gamma}.
$$
where $2<\gamma<3$. Therefore, when the number of iterations $t$ is large enough, the iterated network $G_t$ generated by this weighted $m-$clique annex operation $\tau^r_m(\cdot)$ is a scale-free network model with power exponent $\gamma=2+\frac{\ln(\frac{m}{2}+1)}{\ln(m+1)}$.\cite{barabasi1999emergence}

\subsection{Transfractal dimension of the iterative network $G_t$}

Due to the global self-similarity of the network, the network $G_t$ will have transfinite fractal properties when the number of iterations $t \rightarrow \infty$.\cite{Rozenfeld_2007} First, based on the network architecture, it is easy to prove that the network diameter $L_{t}$ satisfies:
$$
L_{t}=L_{0}+t,
$$
where, $L_{0}$ is the diameter of the initial network $G_0$.
Then, the network transfinite fractal dimension $\bar{d}_{f}$ satisfies \cite{Rozenfeld_2007}:
$$
|V_{L_{t}+\ell}|=e^{\ell \bar{d}_{f}} |V_{L_{t}}|,~~~ t \rightarrow \infty.
$$
Therefore, the transfinite fractal dimension $\bar{d}_{f}$ of the iterative network $G_t$ is:
$$
\bar{d}_{f}=\ln(\frac{1}{2}m^{2}+\frac{3}{2}m+1).
$$

\section{Spectral analysis}

\subsection{Definition of matrix and its spectrum}

In this Section, the spectral properties of the normalized Laplacian matrix corresponding to the iterative weighted network $G_t$ will be studied. Firstly, the relevant definitions need to be given. For the iterative weighted network $G_t$ $(t\geq0)$, the corresponding generalized adjacency matrix, denoted as $W_t$, is the square matrix with $|V_t|$ rows and $|V_t|$ columns, where the element of the $i$-th row and the $j$-th column is denoted as $W_t(i,j)$, which satisfies the following relation:
$$
W_t(i,j)=\left\{
\begin{array}{rcl}
\omega_{ij} ~~& , & i\sim j \\
0 ~~& , & \textrm{otherwise}
\end{array}
\right.,
$$
where $i$ and $j$ correspond to two nodes in the network $G_t$ respectively, while $i\sim j$ represents that two nodes are adjacent. Then, the diagonal strength matrix of iterative weighted network $G_t$ can be defined as: $D_t=diag\{s_t^1,s_t^2,\ldots,s_t^{|V_t|}\}$. Based on the generalized adjacency matrix $W_t$ and diagonal strength matrix $D_t$, the probability transfer matrix, denoted as $T_t$, can be defined as: $T_t=D_t^{-1}W_t$, where the element $T_t(i,j)$ in the $i$-th row and the $j$-th column represents the probability of the walker transferring from node $i$ to node $j$. It is noteworthy that the probability transfer matrix $T_t$ is an asymmetric matrix in most cases. In order to make it symmetric, the normalized adjacency matrix $P_t$ of iterative weighted network $G_t$ can be defined as:
\begin{eqnarray}\label{eq:3.1.0}
P_t=D_t^{-\frac{1}{2}}W_tD_t^{-\frac{1}{2}}=D_t^{-\frac{1}{2}}T_tD_t^{\frac{1}{2}}.
\end{eqnarray}
Clearly, the normalized adjacency matrix $P_t$ is similar to the probability transfer matrix $T_t$, so they have the same eigenvalues. Then, we can define the normalized Laplacian matrix $L_t$ as:
\begin{eqnarray}\label{eq:3.1.1}
L_t=I_{|V_t|}-D_t^{-\frac{1}{2}}T_tD_t^{\frac{1}{2}}=I_{|V_t|}-P_t,
\end{eqnarray}
where $I_{|V_t|}$ is the identity matrix with rank $|V_t|$.

To represent the spectrum of normalized Laplacian matrix of the iterative weighted network $G_t$, the eigenvalues of matrix $L_t$ are denoted as $\sigma_t^i$ $(1\leq i\leq|V_t|)$ which satisfy $\sigma_t^1\leq\sigma_t^2\leq\cdots\leq\sigma_t^{|V_t|}$, and the set of all its eigenvalues is denoted as $\Psi_t=\{\sigma_t^1,\sigma_t^2,\ldots,\sigma_t^{|V_t|}\}$. Similarly, the eigenvalues of normalized adjacency matrix $P_t$ of the network $G_t$ are defined as $\lambda_t^i$ $(1\leq i\leq|V_t|)$, which satisfies $\lambda_t^1\geq\lambda_t^2\geq\cdots\geq\lambda_t^{|V_t|}$, and the set of all eigenvalues is denoted as: $\Lambda_t=\{\lambda_t^1,\lambda_t^2,\ldots,\lambda_t^{|V_t|}\}$. From the Eq.(\ref{eq:3.1.1}), it is easy to prove that its eigenvalues satisfy the following relation:
\begin{eqnarray}\label{eq:3.1.2}
\sigma_t^i=1-\lambda_t^{i}, ~~~ i\in\big[1,|V_t|\big].
\end{eqnarray}
Therefore, only by solving all the eigenvalues of the normalized adjacency matrix $P_t$ of the iterative weighted network $G_t$, the spectrum of normalized Laplacian matrix can be given by Eq.(\ref{eq:3.1.2}). In addition, previous studies have proved that the eigenvalue $\lambda_t^i$ of normalized adjacency matrix $P_t$ of weighted connected network must satisfy the following properties:\cite{chung1997spectral}
$$
1=\lambda_t^1>\lambda_t^2\geq\cdots\geq\lambda_t^{|V_t|}\geq-1.
$$
Therefore, the spectrum of the normalized Laplacian matrix $L_t$ of the iterative weighted network $G_t$ that satisfies the following condition:
$$
0=\sigma_t^1<\sigma_t^2\leq\cdots\leq\sigma_t^{|V_t|}\leq2,
$$
whih is very important for us to solve the eigenvalues of the matrices $P_t$ and $L_t$ later.

\subsection{Iteration relation of spectrum of the Laplacian matrix}

Since the iterative weighted network $G_t$ is generated by the operation $\tau^r_m(\cdot)$, the spectrum of the normalized Laplacian matrix between two successive generations may also have some correlation. In this Subsection, the iterative relationship of the eigenvalues of the normalized Laplacian matrix $L_t$ will be revealed.

For any eigenvalue $\lambda_t\in\Lambda_t$ of the normalized adjacency matrix $P_t$, the corresponding eigenvector is denoted as $\psi_t$, where $\psi_t$ can be expressed as $\psi_t=[\psi_t^1,\psi_t^2,\ldots,\psi_t^{|V_t|}]^\top$, and the element $\psi_t^i$ in the vector corresponds to the node $i$ in the network $G_t$. According to the definition given by Eq.(\ref{eq:3.1.0}), the elements in the normalized adjacency matrix can be expressed as:
$$
P_t(i,j)=\frac{\omega_{i,j}}{\sqrt{s_t^i\cdot s_t^j}}.
$$
If node $i$ is specified to satisfy $i\in \hat V_t$, the following equation can be obtained from the relationship between eigenvalues and eigenvectors:
\begin{eqnarray}\label{eq:3.2.1}
\lambda_t\cdot\psi_t^i&=&\sum_{j\in V_t^i}P_t(i,j)\psi_t^j\nonumber\\
&=&\sum_{j\in\hat V_t^i}\frac{\omega_{ij}}{\sqrt{s_t^i\cdot s_t^j}}\psi_t^j+\sum_{k\in \bar V_t^i}\frac{\omega_{i,k}}{\sqrt{s_t^i\cdot s_t^k}}\psi_t^k.
\end{eqnarray}
From Eq(\ref{eq:2.2.3}), the first term of Eq.(\ref{eq:3.2.1}) can be simplified to obtain the following relation:
\begin{eqnarray}\label{eq:3.2.1.1}
\sum_{j\in\hat V_t^i}\frac{\omega_{ij}}{\sqrt{s_t^i\cdot s_t^j}}\psi_t^j=\frac{1}{mr+1} \sum_{j\in V_{t-1}^i}\frac{\omega_{ij}}{\sqrt{s_{t-1}^i\cdot s_{t-1}^j}}\psi_t^j.
\end{eqnarray}
When $k\in\bar V_t^i$, it can be obtained from the definition in the first Section that
$$
\bar V_t^i=\bigcup_{j\in\hat V_t^i}\bar V_t^{ij}.
$$
Therefore, the second term of Eq.(\ref{eq:3.2.1}) can be expressed as:
\begin{eqnarray}\label{eq:3.2.2}
\sum_{k\in \bar V_t^i}\frac{\omega_{ik}}{\sqrt{s_t^i\cdot s_t^k}}\psi_t^k = \sum_{j\in\hat V_t^i}\sum_{k\in \bar V_t^{ij}} \frac{\omega_{ik}}{\sqrt{s_t^i\cdot s_t^k}}\psi_t^k
\end{eqnarray}
Obviously, the number of nodes in $\bar V_t^{ij}$ is $m$ and these nodes are equivalent to each other, so it is easy to prove that Eq.(\ref{eq:3.2.2}) can be simplified to the following form:
\begin{eqnarray}\label{eq:3.2.3}
\sum_{k\in \bar V_t^i}\frac{\omega_{ik}}{\sqrt{s_t^i\cdot s_t^k}}\psi_t^k &=& \sum_{j\in\hat V_t^i} \sum_{k\in \bar V_t^{ij}} \frac{\omega_{ik}}{\sqrt{s_t^i\cdot s_t^k}}\psi_t^k
=\sum_{j\in\hat V_t^i}\frac{mr\omega_{ij}}{\sqrt{s_t^i\cdot (m+1)r\omega_{ij}}}\psi_t^k.
\end{eqnarray}
Then, when $k\in \bar V_t^{ij}$, the relation between its corresponding eigenvector and eigenvalue can be expressed as:
\begin{eqnarray}\label{eq:3.2.4}
\lambda_t\cdot\psi_t^k&=&\sum_{h\in V_t^k}P_t(k,h)\psi_t^h\nonumber\\
&=&\frac{\omega_{ik}}{\sqrt{s_t^i\cdot s_t^k}}\psi_t^i+\frac{\omega_{jk}}{\sqrt{s_t^j\cdot s_t^k}}\psi_t^j
+\sum_{h\in \bar V_t^{ij}/\{k\}}\frac{\omega_{kh}}{\sqrt{s_t^k\cdot s_t^h}}\psi_t^h\nonumber\\
&=&\frac{\omega_{ik}}{\sqrt{s_t^i\cdot s_t^k}}\psi_t^i+\frac{\omega_{jk}}{\sqrt{s_t^j\cdot s_t^k}}\psi_t^j
+\frac{m-1}{m+1}\psi_t^k,
\end{eqnarray}
where we used the conclusions: $\omega_{kh}=r\omega_{ij}$, $s_t^k=s_t^h=(m+1)r\omega_{ij}$ and the equivalence relation of the inner nodes of set $\bar V_t^{ij}$. Therefore, the eigenvector element $\psi_t^k$ corresponding to node $k\in \bar V_t^{ij}$ satisfies the following relation:
\begin{eqnarray}\label{eq:3.2.5}
\psi_t^k=\frac{m+1}{\lambda_t(m+1)-m+1}\cdot\frac{r\omega_{ij}}{\sqrt{s_t^k}}\Big(\frac{\psi_t^i}{\sqrt{s_t^i}} + \frac{\psi_t^j}{\sqrt{s_t^j}}\Big),
\end{eqnarray}
where the eigenvalue $\lambda_t$ must satisfy: $\lambda_t\neq\frac{m-1}{m+1}$. Then, the eigenvector element $\psi_t^k$ in Eq.(\ref{eq:3.2.5}) is substituted into Eq.(\ref{eq:3.2.3}) to obtain:
\begin{eqnarray}\label{eq:3.2.6}
\sum_{k\in \bar V_t^i}\frac{\omega_{ik}}{\sqrt{s_t^i\cdot s_t^k}}\psi_t^k &=&
\sum_{j\in\hat V_t^i} \frac{mr\omega_{ij}}{\lambda_t(m+1)-m+1}\Big(\frac{\psi_t^i}{s_t^i} + \frac{\psi_t^j}{\sqrt{s_t^is_t^j}}\Big)\nonumber\\
&=&\frac{1}{mr+1}\sum_{j\in V_{t-1}^i} \frac{mr\omega_{ij}}{\lambda_t(m+1)-m+1}\Big(\frac{\psi_t^i}{s_{t-1}^i} + \frac{\psi_t^j}{\sqrt{s_{t-1}^is_{t-1}^j}}\Big).~~~
\end{eqnarray}
Therefore, if Eq.(\ref{eq:3.2.1.1}) and Eq.(\ref{eq:3.2.6}) are substituted into Eq.(\ref{eq:3.2.1}), the following relation can be obtained:
\begin{eqnarray}\label{eq:3.2.7}
\lambda_t\cdot\psi_t^i&=&\frac{1}{mr+1} \sum_{j\in V_{t-1}^i} \Big[ \frac{mr\omega_{ij}}{\lambda_t(m+1)-m+1}\Big(\frac{\psi_t^i}{s_{t-1}^i} + \frac{\psi_t^j}{\sqrt{s_{t-1}^is_{t-1}^j}}\Big)+\frac{\omega_{ij}}{\sqrt{s_{t-1}^i\cdot s_{t-1}^j}}\psi_t^j\Big]\nonumber\\
&=&\frac{\lambda_t(m+1)-m+1+mr}{(\lambda_t(m+1)-m+1)(mr+1)}\sum_{j\in V_{t-1}^i}\frac{\omega_{ij}}{\sqrt{s_{t-1}^i\cdot s_{t-1}^j}}\psi_t^j +\frac{mr\psi_t^j}{(\lambda_t(m+1)-m+1)(mr+1)}.
\end{eqnarray}
Then simplify Eq.(\ref{eq:3.2.7}) to obtain the following equation:
\begin{eqnarray}\label{eq:3.2.8}
\frac{\lambda_t(\lambda_t(m+1)-m+1)(mr+1)-mr}{\lambda_t(m+1)-m+1+mr}\psi_t^i=\sum_{j\in V_{t-1}^i}P_{t-1}(i,j)\psi_t^j,
\end{eqnarray}
where $\lambda_t\neq\frac{m-1}{m+1}$ and $\lambda_t\neq\frac{m-1-mr}{m+1}$.
From Eq.(\ref{eq:3.2.8}), it is easy to prove that when $i\in\hat V_t$, the new vector $\psi_t^*$ composed of the eigenvector element $\psi_t^i$ satisfies the following relation:
\begin{eqnarray}\label{eq:3.2.9}
\frac{\lambda_t(\lambda_t(m+1)-m+1)(mr+1)-mr}{\lambda_t(m+1)-m+1+mr}\psi_t^*\nonumber\triangleq\lambda_t^*\cdot\psi_t^*=P_{t-1}\cdot\psi_t^*.
\end{eqnarray}
Thus, it can be determined that $\lambda_t^*$ is the eigenvalue of normalized adjacency matrix $P_{t-1}$, and $\psi_t^*$ is the eigenvector corresponding to the eigenvalue $\lambda_t^*$. The above Eq.(\ref{eq:3.2.9}) reflects the relationship between the eigenvalues of two successive generations of normalized adjacency matrices.
Thus, it is easy to prove that, for $\forall \lambda_{t-1}\in\Lambda_{t-1}$, there is $\lambda_t\in\Lambda_t$ such that the following relationship holds:
\begin{eqnarray}\label{eq:3.2.10}
(mr+1)(m+1)\lambda_t^2-\big[(m-1)(mr+1)+(m+1)\lambda_{t-1}\big]\lambda_t
-\big[mr+(mr+1-m)\lambda_{t-1}\big]=0
\end{eqnarray}
Then, we define the equations $f_1(x)$ and $f_2(x)$ as follows:
\begin{align}
f_1(x)=&\frac{1}{2(mr+1)(m+1)}\Bigg[(m-1)(mr+1)+(m+1)x+\nonumber\\
&\sqrt{\big[(m-1)(mr+1)+(m+1)x\big]^2+4(mr+1)(m+1)\big[mr+(mr+1-m)x\big]}\Bigg],\label{eq:3.2.11}\\
f_2(x)=&\frac{1}{2(mr+1)(m+1)}\Bigg[(m-1)(mr+1)+(m+1)x-\nonumber\\
&\sqrt{\big[(m-1)(mr+1)+(m+1)x\big]^2+4(mr+1)(m+1)\big[mr+(mr+1-m)x\big]}\Bigg].\label{eq:3.2.12}
\end{align}
Therefore, it can be concluded that for $\forall \lambda_{t-1}\in\Lambda_{t-1}$, $f_1(\lambda_{t-1})$ and $f_2(\lambda_{t-1})$ are the eigenvalues of the adjacency matrix $P_t$, that is $f_1(\lambda_{t-1}),f_2(\lambda_{t-1})\in\Lambda_t$.
Here, it is worth noting that, according to the properties of matrix eigenvalues in Section 3.1 above, there is only one maximum eigenvalue in $\Lambda_{t-1}$ that satisfies: $\lambda_{t-1}^1=1$, so it can be obtained that:
$$
f_1(\lambda_{t-1}^1)=1 ~~~\textrm{and}~~~ f_2(\lambda_{t-1}^1)=-\frac{2mr+1-m}{(mr+1)(m+1)}.
$$
The monotonicity of $f_1(x)$ and $f_2(x)$ also ensures that the eigenvalues generated by $f_1(x)$ and $f_2(x)$ meet the following requirements:
$$
1>f_1(\lambda_{t-1})>f_2(\lambda_{t-1})\geq-1 ~~~\textrm{where}~~~ \lambda_{t-1}\in\Lambda_{t-1}/\{\lambda_{t-1}^1\}.
$$

Then, combining Eq(\ref{eq:3.1.2}) and Eq(\ref{eq:3.2.10}), it is easy to prove that the eigenvalues of the Laplacian matrix $L_t$ satisfy the following relationship:
\begin{eqnarray}\label{eq:3.2.13}
(mr+1)(m+1)\sigma_t^2-\big[(m+1)\sigma_{t-1}+mr(m+3)+2\big]\sigma_t+(m+2)\sigma_{t-1}=0.
\end{eqnarray}
Similarly, equations $g_1(x)$ and $g_2(x)$ can be defined as follows:
\begin{eqnarray}
g_1(x)&=&\frac{1}{2(mr+1)(m+1)}\Bigg[(m+1)\sigma_{t-1}+mr(m+3)+2\nonumber\\
&&+\sqrt{\big[(m+1)\sigma_{t-1}+mr(m+3)+2\big]^2-4(mr+1)(m+1)(m+2)\sigma_{t-1}}\Bigg],\label{eq:3.2.14}\\
g_2(x)&=&\frac{1}{2(mr+1)(m+1)}\Bigg[(m+1)\sigma_{t-1}+mr(m+3)+2\nonumber\\
&&-\sqrt{\big[(m+1)\sigma_{t-1}+mr(m+3)+2\big]^2-4(mr+1)(m+1)(m+2)\sigma_{t-1}}\Bigg].\label{eq:3.2.15}
\end{eqnarray}
Therefore, for $\forall\sigma_{t-1}\in\Psi_{t-1}$, $g_1(\sigma_{t-1})$ and $g_2(\sigma_{t-1})$ are the eigenvalues of the normalized Laplacian matrix $L_t$ of the $t$ generation respectively, that is $g_1(\sigma_{t-1}),g_2(\sigma_{t-1})\in\Psi_t$. In addition, for $\sigma_{t-1}^1=1$, it can be obtained that
$$
g_1(\sigma_{t-1}^1)=\frac{m^2r+3mr+2}{(mr+1)(m+1)} ~~~\textrm{and}~~~ g_2(\sigma_{t-1}^1)=0.
$$
For $\sigma_{t-1}\in\Psi_{t-1}/\{\sigma_{t-1}^1\}$, the following relation can be obtained from the monotonicity of $g_1(x)$ and $g_2(x)$:
$$
0<g_2(\sigma_{t-1})<g_1(\sigma_{t-1})\leq 2.
$$

Hence, iterative relation of eigenvalues of the normalized Laplacian matrix $L_t$ has been derived. But it's worth noting that Eq.(\ref{eq:3.2.14}) and Eq.(\ref{eq:3.2.15}) indicate that each eigenvalue in the matrix $L_{t-1}$ will correspond to two eigenvalues in the matrix $L_t$. Therefore, the iterative relationship can only determine $2|V_{t-1}|$ eigenvalues in the matrix $L_t$, and there are still $|V_t|-2|V_{t-1}|$ eigenvalues need to be confirmed.

\subsection{Spectrum of Laplacian matrix of weighted iterative network $G_t$}

According to the definition of weighted $m-$clique annex operation $\tau^r_m(\cdot)$, it is easy to prove that the normalized adjacency matrix $P_t$ corresponding to the weighted iterative network $G_t$ satisfies the following iterative relation:
\begin{equation}\label{eq:3.3.1}
P_t=\left[
  \begin{array}{ccccc}
    \frac{P_{t-1}}{mr+1} & X_1 & X_2 & \cdots & X_{|E_{t-1}|} \\
    X_1^\top & A & 0 & \cdots & 0\\
    X_2^\top & 0 & A & \cdots & 0\\
    \vdots & \vdots & \vdots & \ddots & \vdots\\
    X_{|E_{t-1}|}^\top & 0 & 0 & \cdots & A
  \end{array}
\right]
\end{equation}
where the matrix $X_e (1\leq e \leq |E_{t-1}|)$ is a matrix with $|V_{t-1}|$ rows and $m$ columns, and $A$ is the square matrix with $m$ rows. In the matrix $X_e$, $e$ corresponds to an edge in the set $E_{t-1}$, and the two nodes connected by this edge are denoted as $i_e$ and $j_e$, respectively. Therefore, the all elements of row $i_e$ and the all elements of row $j_e$, denoted as $X_e(i_e,k)$ and $X_e(j_e,k)$ $1\leq k \leq m$ respectively, of matrix $X_e$ satisfy:
\begin{eqnarray*}
X_e(i_e,k)=\sqrt{\frac{r\omega_{ij}}{(mr+1)(m+1)s^{i_e}_{t-1}}},\\
X_e(j_e,k)=\sqrt{\frac{r\omega_{ij}}{(mr+1)(m+1)s^{j_e}_{t-1}}}
\end{eqnarray*}
and the other elements are all $0$. In matrix $A$, the diagonal elements are all $0$, and the other elements are all $\frac{1}{m+1}$. It is easy to prove that the eigenvalues of matrix $A$ are $-\frac{1}{m+1}$ and $\frac{m-1}{m+1}$, respectively, and their multiples are $m$ and $1$. Since matrix $A$ reflects the local structure characteristics of network $G_t$, we assume that the eigenvalue of matrix $A$ is also the eigenvalue of normalized adjacency matrix $P_t$. Through the relationship between the matrix and the eigenvalue, it can be proved that if the $t-1$th generation network $G_{t-1}$ satisfies $|E_{t-1}|\geq|V_{t-1}|$, both $-\frac{1}{m+1}$ and $\frac{m-1}{m+1}$ are the eigenvalues of the matrix $P_t$, and their multiples are $|V_t|-|E_{t-1}|-|V_{t-1}|$ and $|E_{t-1}|-|V_{t-1}|$ respectively; if the $t-1$th generation network $G_{t-1}$ satisfies $|E_{t-1}|<|V_{t-1}|$, $-\frac{1}{m+1}$ is the eigenvalue of the matrix $P_t$, and its multiple is $|V_t|-2|V_{t-1}|$. Obviously, according to the fundamental theorem of connected network, if $|E_{t-1}|<|V_{t-1}|$, there must be $|V_{t-1}|=|E_{t-1}|+1$, in terms of network structure, which means that there is no loop in network $G_{t-1}$, that is, $G_{t-1}$ is a tree. Then, considering the definition of network operation $\tau^r_m(\cdot)$, there must be a loop structure on the network after one iteration. Therefore, only the initial network $G_0$ may be a tree, so it needs to be divided into two categories for discussion.

\subsubsection{$G_0$ is not a tree or $t>1$}

Based on the above analysis, we can divide the set $\Lambda_{t}$ composed of the eigenvalues of the matrix $P_t$ into four subsets, which can be expressed as:
\begin{eqnarray}\label{eq:3.3.1.1}
\Lambda_{t}=\Lambda_{t}^1\cup\Lambda_{t}^2\cup\Lambda_{t}^3\cup\Lambda_{t}^4,
\end{eqnarray}
where
\begin{eqnarray*}
\Lambda_{t}^1&=&\{1,-\frac{2mr+1-m}{(mr+1)(m+1)}\},~~
\Lambda_{t}^2=\{f_1(\lambda_{t-1}^2),f_2(\lambda_{t-1}^2),\ldots,f_1(\lambda_{t-1}^{|V_{t-1}|}),f_2(\lambda_{t-1}^{|V_{t-1}|})\},\\
\Lambda_{t}^3&=&\underbrace{\{-\frac{1}{m+1},\ldots,-\frac{1}{m+1}\}}_{|V_t|-|V_{t-1}|-|E_{t-1}|},~~\textrm{and}~~
\Lambda_{t}^4=\underbrace{\{\frac{m-1}{m+1},\ldots,\frac{m-1}{m+1}\}}_{|E_{t-1}|-|V_{t-1}|}.
\end{eqnarray*}

Thus, the spectrum $\Psi_t$ of the normalized Laplacian matrix $L_t$ can also be obtained and expressed as:
\begin{eqnarray}\label{eq:3.3.1.2}
\Psi_{t}=\Psi_{t}^1\cup\Psi_{t}^2\cup\Psi_{t}^3\cup\Psi_{t}^4,
\end{eqnarray}
where
\begin{eqnarray*}
\Psi_{t}^1&=&\{0,\frac{m^2r+3mr+2}{(mr+1)(m+1)}\},~~
\Psi_{t}^2=\{g_1(\sigma_{t-1}^2),g_2(\sigma_{t-1}^2),\ldots,g_1(\sigma_{t-1}^{|V_{t-1}|}),g_2(\sigma_{t-1}^{|V_{t-1}|})\},\\
\Psi_{t}^3&=&\underbrace{\{\frac{m+2}{m+1},\ldots,\frac{m+2}{m+1}\}}_{|V_t|-|V_{t-1}|-|E_{t-1}|},~~\textrm{and}~~
\Psi_{t}^4=\underbrace{\{\frac{2}{m+1},\ldots,\frac{2}{m+1}\}}_{|E_{t-1}|-|V_{t-1}|}.
\end{eqnarray*}

\subsubsection{$G_0$ is a tree and $t=1$}

According to the previous analysis, in network $G_1$ generated by tree network $G_0$ after one operation $\tau_m^r(\cdot)$, the spectrum $\Lambda_1$ of its normalized adjacency matrix $P_1$ can be expressed as:
\begin{eqnarray}\label{eq:3.3.2.1}
\Lambda_1=\Lambda_1^1\cup\Lambda_1^2\cup\Lambda_1^3,
\end{eqnarray}
where
\begin{eqnarray*}
\Lambda_1^1&=&\{1,-\frac{2mr+1-m}{(mr+1)(m+1)}\},~~
\Lambda_1^2=\{f_1(\lambda_0^2),f_2(\lambda_0^2),\ldots,f_1(\lambda_0^{|V_0|}),f_2(\lambda_0^{|V_0|})\},\\
~~\textrm{and}~~\Lambda_1^3&=&\underbrace{\{-\frac{1}{m+1},\ldots,-\frac{1}{m+1}\}}_{|V_1|-2|V_0|}.
\end{eqnarray*}

Therefore, under the above conditions, the set of eigenvalues of normalized Laplacian matrix $L_1$ can be divided into three subsets, which can be expressed as:
\begin{eqnarray}\label{eq:3.3.2.2}
\Psi_1=\Psi_1^1\cup\Psi_1^2\cup\Psi_1^3,
\end{eqnarray}
where
\begin{eqnarray*}
\Psi_1^1&=&\{0,\frac{m^2r+3mr+2}{(mr+1)(m+1)}\},~~
\Psi_1^2=\{g_1(\sigma_0^2),g_2(\sigma_0^2),\ldots,g_1(\sigma_0^{|V_0|}),g_2(\sigma_0^{|V_0|})\},\\
~~\textrm{and}~~\Psi_1^3&=&\underbrace{\{\frac{m+2}{m+1},\ldots,\frac{m+2}{m+1}\}}_{|V_1|-2|V_0|}.
\end{eqnarray*}

Hence, when the spectrum $\Psi_0$ of the initial network $G_0$ is known, the spectrum $\Psi_t$ of the iterative weighted network $G_t$ can be obtained from the above relation.

\begin{remark}
The spectrum of the dynamic growth network $\tau^{r}_{m}(G_0)$ shows that: (1) The eigenvalues corresponding to the original network $G_0$ will be controlled by the change rule $g_1(\cdot)$ and $g_2(\cdot)$ after one iteration, which are completely determined by the weight factor $r$ and the scale factor $m$; (2) The remaining eigenvalues are completely determined by the scale factor $m$. Based on the importance of the eigenvalues of the Laplacian matrix, we can infer that the local structure $K_{m}$ introduced here will have a profound and important impact on the overall characteristics of $\tau^{r}_{m}(G_0)$. Furthermore, it can be deduced that in the growth process of the actual dynamic growth network, the local topology of the newly added node block will have an important impact on the overall properties.
\end{remark}

\section{Application of the spectrum}

After obtaining the spectrum of the iterated weighted network $G_t$, the related applications will be presented in this Section, including the Kemeny constant $K(G_t)$, the  multiplicative degree-Kirchhoff index $\tilde{\mathcal{K}}(G_t)$, and the number of weighted spanning trees $N_t$.

\subsection{Kemeny constant $K(G_t)$}

The random walk on the weighted network $G_t$ is defined as: if the walker starts from node $i$, then after a unit time, the walker will transfer from node $i$ to its neighbor node $j\in V_t^i$ with probability $\frac{\omega_{ij}}{s_t^i}$. Reviewing the definition of probability transfer matrix $T_t$ of network $G_t$, it can be found that the random walk process can be completely described by this matrix, where the element $T_t(i,j)$ in the $i$-th row and $j$-th column of the matrix represents the probability of the particle moving from node $i$ to node $j$.

Then, the mean first-passage time $F_{ij}$ from the node $i$ to the node $j$ on the weighted network $G_t$ can be defined as: the expectation of the time required for the random walker departing from node $i$ and arriving at node $j$ for the first time.\cite{redner2001guide} In addition, the stationary distribution $\pi_t$ on the weighted network $G_t$ be defined as:\cite{lovasz1993random,aldous1995reversible}
$$
\pi_t=[\pi_t^1,\pi_t^2,\ldots,\pi_t^{|V_t|}],
$$
where for $1\leq i\leq|V_t|$, $\pi_t^i=\frac{s_t^i}{2S_t}$. It is easy to prove that the stationary distribution $\pi_t$ satisfies the following two properties:
$$
\pi_t T_{t}=\pi_t ~~~\textrm{and}~~~ \sum_{i=1}^{|V_t|}\pi_t^i=1.
$$

Thus, Kemeny constant $K(G_t)$ on the weighted network $G_t$ can be defined as the mean value of the first passage time $F_{ij}$ from any node $i$ to the target node $j$, where the target node $j$ is randomly selected according to the stationary distribution $\pi_t$. The mathematical expression of this definition is:\cite{aldous1995reversible}
$$
K(G_t)=\sum_{j=1}^{|V_t|}\pi_t^jF_{ij}.
$$
In fact, the Kemeny constant $K(G_t)$ is independent of the selection of starting node $i$, so the constant reflects the global property of random walk on weighted network $G_t$. Moreover, previous work has proved that Kemeny constant $K(G_t)$ can be completely determined by the reciprocal sum of non-zero eigenvalues of the normalized Laplacian matrix of the network $G_t$, that is:\cite{aldous1995reversible,levene2002kemeny}
\begin{eqnarray}\label{eq:4.1.1}
K(G_t)=\sum_{i=2}^{|V_t|}\frac{1}{1-\lambda_t^i}=\sum_{i=2}^{|V_t|}\frac{1}{\sigma_t^i}.
\end{eqnarray}
Therefore, according to the spectrum analysis of Laplacian matrix of network $G_t$, the analytical expression of $K(G_t)$can be derived.

When the number of iterations $t>1$, the following relationship can be obtained according to the classification of the spectrum in Eq.(\ref{eq:3.3.1.2}):
\begin{eqnarray}\label{eq:4.1.2}
K(G_t)=\frac{(mr+1)(m+1)}{m^2r+3mr+2} + \sum_{\sigma_t\in\Psi_t^2}\frac{1}{\sigma_t} + \sum_{\sigma_t\in\Psi_t^3}\frac{1}{\sigma_t} + \sum_{\sigma_t\in\Psi_t^4}\frac{1}{\sigma_t}.
\end{eqnarray}
Using Vieta formulas for Eq.(\ref{eq:3.2.13}) is easy to prove the following relationship:
\begin{eqnarray}
g_1(\sigma_{t-1})+g_2(\sigma_{t-1})&=&\frac{(m+1)\sigma_{t-1}+mr(m+3)+2}{(mr+1)(m+1)},\label{eq:4.1.3}\\
g_1(\sigma_{t-1})\cdot g_2(\sigma_{t-1})&=&\frac{(mr+2)\sigma_{t-1}}{(mr+1)(m+1)},\label{eq:4.1.4}
\end{eqnarray}
where $\sigma_{t-1}\in \Psi_{t-1}/\{0\}$. Therefore, the sum of the reciprocal of the elements in set $\Psi_t^2$ satisfies the following relation:
\begin{eqnarray}\label{eq:4.1.5}
\sum_{\sigma_t\in\Psi_t^2}\frac{1}{\sigma_t}&=&\sum_{i=2}^{V_{t-1}}\big(\frac{1}{g_1(\sigma_{t-1}^i)} + \frac{1}{g_2(\sigma_{t-1}^i)}\big)=\sum_{i=2}^{V_{t-1}}\Big(\frac{m+1}{mr+2}+\frac{(m+3)mr+2}{mr+2}\cdot\frac{1}{\sigma_{t-1}^i}\Big)\nonumber\\
&=&\frac{2|E_0|(m\!+\!1)}{(mr\!+\!2)(m\!+\!3)}(\frac{m^2\!+\!3m\!+\!2}{2})^{t-1}\!+\!\frac{(m+3)mr+2}{mr+2}\!\cdot\! K(G_{t-1})\!+\!\frac{m+1}{mr+2}(|V_0|\!-\!\frac{2|E_0|}{m+2}\!-\!1).
\end{eqnarray}

Similarly, through the spectral analysis in Eq.(\ref{eq:3.3.1.2}), the third and fourth terms in Eq.(\ref{eq:4.1.2}) can be expressed as:
\begin{eqnarray}
\sum_{\sigma_t\in\Psi_t^3}\frac{1}{\sigma_t}&=&(|V_t|-|E_{t-1}|-|V_{t-1}|)\frac{m+1}{m+2}=|E_0|\frac{m^2-1}{m+2}
(\frac{m^2+3m+2}{2})^{t-1}\label{eq:4.1.6}\\
\sum_{\sigma_t\in\Psi_t^4}\frac{1}{\sigma_t}&=&(|E_{t-1}|-|V_{t-1}|)\frac{m+1}{2}=|E_0|\frac{(m+1)^2}{2(m+3)}(\frac{m^2+3m+2}{2})^{t-1}-\frac{m+1}{2}(|V_0|-\frac{2|E_0|}{m+3}).~~~\label{eq:4.1.7}
\end{eqnarray}
Then, after the Eq.(\ref{eq:4.1.5}), Eq.(\ref{eq:4.1.6}) and Eq.(\ref{eq:4.1.7}) are brought into Eq.(\ref{eq:4.1.2}), the iterative expression of Kemeny constant $K(G_t)$ can be expressed as:
\begin{eqnarray}\label{eq:4.1.8}
K(G_t)&=&\frac{(m+3)mr+2}{mr+2}K(G_{t-1})
+|E_0|\Big(\frac{2(m+1)}{(mr+2)(m+3)}+\frac{m^2-1}{m+2}+\frac{(m+1)^2}{2(m+3)}\Big)\cdot\nonumber\\
&&(\frac{m^2+3m+2}{2})^{t-1}\!+\!\frac{m+1}{mr+2}(|V_0|\!-\!\frac{2|E_0|}{m+2}-1) \!-\! \frac{m+1}{2}(|V_0|-\frac{2|E_0|}{m+3}) \!+\!\frac{(mr+1)(m+1)}{m^2r+3mr+2}.
\end{eqnarray}
From the above iterative relation of Kemeny constant $K(G_t)$, it is easy to deduce the relation between $K(G_t)$ and $K(G_1)$ as follows:
\begin{eqnarray}\label{eq:4.1.9}
K(G_t)&=&\big(\frac{(m+3)mr+2}{mr+2}\big)^{t-1}K(G_1)
+|E_0|\Big(\frac{2(m+1)}{(mr+2)(m+3)}+\frac{m^2-1}{m+2}+\frac{(m+1)^2}{2(m+3)}\Big)\cdot\nonumber\\
&&\sum_{i=0}^{t-2}\big(\frac{(m+3)mr+2}{mr+2}\big)^i\big(\frac{m^2+3m+2}{2}\big)^{t-1-i} +\Big[\frac{m+1}{mr+2}(|V_0|-\frac{2|E_0|}{m+2}-1) \nonumber\\
&&-\frac{m+1}{2}(|V_0|-\frac{2|E_0|}{m+3}) +\frac{(mr+1)(m+1)}{m^2r+3mr+2}\Big]
\sum_{i=0}^{t-2}\big(\frac{(m+3)mr+2}{mr+2}\big)^i.
\end{eqnarray}
Then, from the spectrum analysis in the previous Section, it can be seen that the spectrum of network $G_1$ will be different according to whether $G_0$ is a tree, so it needs to be classified and discussed.

\subsubsection{Case 1: $G_0$ is not a tree}

When the initial network $G_0$ is not a tree, the following relationship between $K(G_1)$ and $K(G_0)$ can be naturally obatined:
\begin{eqnarray}\label{eq:4.1.1.1}
K(G_1)&=&\frac{(m+3)mr+2}{mr+2}K(G_0)
+|E_0|\Big(\frac{2(m+1)}{(mr+2)(m+3)}+\frac{m^2-1}{m+2}+\frac{(m+1)^2}{2(m+3)}\Big)\nonumber\\
&&+\frac{m+1}{mr+2}(|V_0|-\frac{2|E_0|}{m+2}-1)
-\frac{m+1}{2}(|V_0|-\frac{2|E_0|}{m+3}) +\frac{(mr+1)(m+1)}{m^2r+3mr+2}.
\end{eqnarray}
Obviously, by combining Eq.(\ref{eq:4.1.9}) with Eq.(\ref{eq:4.1.1.1}), the relation between $K(G_t)$ and $K(G_0)$ can be deduced, but it is not an analytic expression for $K(G_t)$. To find an analytic expression for $K(G_t)$, the relationship between $\frac{(m+3)mr+2}{mr+2}$ and $\frac{m^2+3m+2}{2}$ must be discussed.
Since scale factor $m\in\mathbb{N}^+$ and weight factor $r\in\mathbb{R}^+$, it is easy to prove the following relation:
\begin{eqnarray}\label{eq:4.1.1.2}
\left\{
\begin{array}{ll}
\frac{(m+3)mr+2}{mr+2}>\frac{m^2+3m+2}{2} , & m=1 ~~\textrm{and}~~ r>4 \\
\frac{(m+3)mr+2}{mr+2}=\frac{m^2+3m+2}{2} , & m=1 ~~\textrm{and}~~ r=4 \\
\frac{(m+3)mr+2}{mr+2}<\frac{m^2+3m+2}{2} , & \textrm{others.}
\end{array}
\right.
\end{eqnarray}
Therefore, when $m\neq1$ or $r\neq4$, we can deduce the analytic expression of Kemeny constant $K(G_t)$ as follows:
\begin{eqnarray}\label{eq:4.1.1.3}
K(G_t)=\Big(K(G_0)-k_1-k_2\Big)\big(\frac{(m+3)mr+2}{mr+2}\big)^t
+k_1\big(\frac{m^2+3m+2}{2}\big)^t+k_2,
\end{eqnarray}
where the coefficients $k_1$ and $k_2$ are determined by the number of nodes $|V_0|$, the number of edges $|E_0|$, the size factor $m$, and the weight factor $r$ in the initial network, which are expressed as:
\begin{eqnarray}
k_1&=&\frac{|E_0|(m+1)(20m-4mr+7m^2r+3m^3r+6m^2+4)}{m(m^2+5m+6)(2m-4r+mr+m^2r+6)},\label{eq:4.1.1.4}\\
k_2&=&\frac{(mr+2)(m+1)}{2(m+2)(m+3)(m^3r^2+3m^2r^2+2m^2r+8mr+4)}
(6|V_0|-4|E_0|+6m+2|V_0|m\nonumber\\
&&-6mr-2m^2r+2m^2-6|E_0|mr+9|V_0|mr-2|E_0|m^2r+6|V_0|m^2r+|V_0|m^3r).\label{eq:4.1.1.5}
\end{eqnarray}

When $m=1$ and $r=4$, it is easy to calculate that the analytic expression of $K(G_t)$ can be expressed as follows:
\begin{eqnarray}\label{eq:4.1.1.6}
K(G_t)=(K(G_0)+\frac{1}{6}|E_0|-\frac{1}{3}|V_0|+\frac{1}{9})\cdot3^t+\frac{|E_0|}{4}3^t\cdot t -(\frac{1}{6}|E_0|-\frac{1}{3}|V_0|+\frac{1}{9}).
\end{eqnarray}
After obtaining the analytic expression of $K(G_t)$, we can naturally obtain that when the number of iteration $t$ is large enough, the main term of Kemeny constant $K(G_t)$ obeys:
\begin{eqnarray}\label{eq:4.1.1.7}
K(G_t)\sim\left\{
\begin{array}{ll}
\Big(K(G_0)-k_1-k_2\Big)\big(\frac{(m+3)mr+2}{mr+2}\big)^t\sim &\\
\Big(K(G_0)-k_1-k_2\Big)\Big(\frac{(m+3)}{2|E_0|}|V_t|\Big)^{\ln{ \frac{2(m+3)mr+4}{(mr+2)(m^2+3m+2)}}},
& m=1 ~~\textrm{and}~~ r>4 \\
\frac{|E_0|}{4}3^t\cdot t \sim \frac{|V_t|\ln|V_t|}{2\ln3}, & m=1 ~~\textrm{and}~~ r=4 \\
k_1\big(\frac{m^2+3m+2}{2}\big)^t\sim\frac{k_1(m+3)}{2|E_0|}|V_t|, & \textrm{others.}
\end{array}
\right.
\end{eqnarray}

Therefore, when the initial network $G_0$ is not a tree, the analytical expression and approximate expression of Kemeny constant $K(G_t)$ on the weighted iterative network $G_t$ have been fully figured out.

\subsubsection{Case 2: $G_0$ is a tree}

When the initial network $G_0$ is a tree, it is easy to prove the relationship between $K(G_1)$ and $K(G_0)$ by the spectral relationship in Eq.(\ref{eq:3.3.2.2}) and $E_0=N_0-1$ as follows:
\begin{eqnarray}\label{eq:4.1.2.1}
K(G_1)&=&(|V_1|-2|V_0|)\frac{m+1}{m+2}+\sum_{\sigma_0\in\Psi_0/\{0\}}\Big(\frac{m+1}{mr+2}+\frac{m^2r+3mr+2}{(mr+2)\sigma_0}\sigma_0\Big) +\frac{(mr+1)(m+1)}{m^2r+3mr+2}\nonumber\\
&=&\frac{m^2r+3mr+2}{mr+2}K(G_0)\!+\!\Big[(m-1)|V_0|-m\Big]\frac{m+1}{m+2}\!+\!(|V_0|-1)\frac{m+1}{mr+2}\!+\!\frac{(mr+1)(m+1)}{m^2r+3mr+2}.
\end{eqnarray}

Similarly, when calculating the analytic expression of Kemeny constant $K(G_t)$, the magnitude relation between $\frac{(m+3)mr+2}{mr+2}$ and $\frac{m^2+3m+2}{2}$ in Eq.(\ref{eq:4.1.1.2}) need to be considered. Therefore, when $m\neq1$ or $r\neq4$, from Eq.(\ref{eq:4.1.9}) and $|E_0|=|V_0|-1$, $K(G_t)$ can be expressed as:
\begin{eqnarray}\label{eq:4.1.2.2}
K(G_t)=\Big(K(G_0)+k_5\Big)\big(\frac{(m+3)mr+2}{mr+2}\big)^t
+k_3\big(\frac{m^2+3m+2}{2}\big)^{t}+k_4,
\end{eqnarray}
where the coefficients $k_3$ and $k_4$ can be expressed as:
\begin{eqnarray}
k_3&=&\frac{(|V_0|-1)(m+1)(20m-4mr+7m^2r+3m^3r+6m^2+4)}{m(m^2+5m+6)(2m-4r+mr+m^2r+6)},\label{eq:4.1.2.3}\\
k_4&=&\frac{(mr+2)(m+1)^2(2|V_0|+2m+3|V_0|mr+|V_0|m^2r+4)}{2(m+2)(m+3)(m^3r^2+3m^2r^2+2m^2r+8mr+4)},\label{eq:4.1.2.4}\\
k_5&=&\frac{mr+2}{(m+3)mr+2}\Big(\Big[(m-1)|V_0|-m\Big]\frac{m+1}{m+2}+(|V_0|-1)\frac{m+1}{mr+2}+\frac{(mr+1)(m+1)}{m^2r+3mr+2}-k_3-k_4\Big).
\end{eqnarray}
In addition, in the special case, namely $m=1$ and $r=4$, it is easy to prove that the analytic expression of $K(G_t)$ is:
\begin{eqnarray}\label{eq:4.1.2.5}
K(G_t)=\frac{|V_0|-1}{4}3^t\cdot t+(K(G_0)-\frac{29}{144}|V_0|+\frac{35}{432})3^t+\frac{3|V_0|+1}{16}.
\end{eqnarray}
Hence, when the initial network $G_0$ is a tree, the analytic expression of Kemeny constant $K(G_t)$ on weighted iterative network $G_t$ have been solved. Similarly, according to the analytic expression, when the number of iterations $t$ is sufficiently large, the main term of $K(G_t)$ obeys the following approximate relation:
\begin{eqnarray}\label{eq:4.1.2.6}
K(G_t)\sim\left\{
\begin{array}{lr}
\Big(K(G_0)+k_5\Big)\big(\frac{(m+3)mr+2}{mr+2}\big)^t\sim &\\
\Big(K(G_0)+k_5\Big)\Big(\frac{(m+3)}{2|E_0|}|V_t|\Big)^{\ln{ \frac{2(m+3)mr+4}{(mr+2)(m^2+3m+2)}}} ,
& m=1 ~~\textrm{and}~~ r>4 \\
\frac{|V_0|-1}{4}3^t\cdot t \sim \frac{|V_t|\ln|V_t|}{2\ln3} , & m=1 ~~\textrm{and}~~ r=4 \\
k_3\big(\frac{m^2+3m+2}{2}\big)^t\sim\frac{k_3(m+3)}{2|E_0|}|V_t| , & \textrm{others.}
\end{array}
\right.
\end{eqnarray}
Therefore, for arbitrary initial network $G_0$, scale factor $m\in\mathbb{N}^+$, weight factor $r\in\mathbb{R}^+$, the analytic expression and approximate expression of Kemeny constant $K(G_t)$ have been fully derived.

\begin{figure}[t]
\centering
\begin{minipage}[t]{0.48\textwidth}
\centering
\includegraphics[width=8cm]{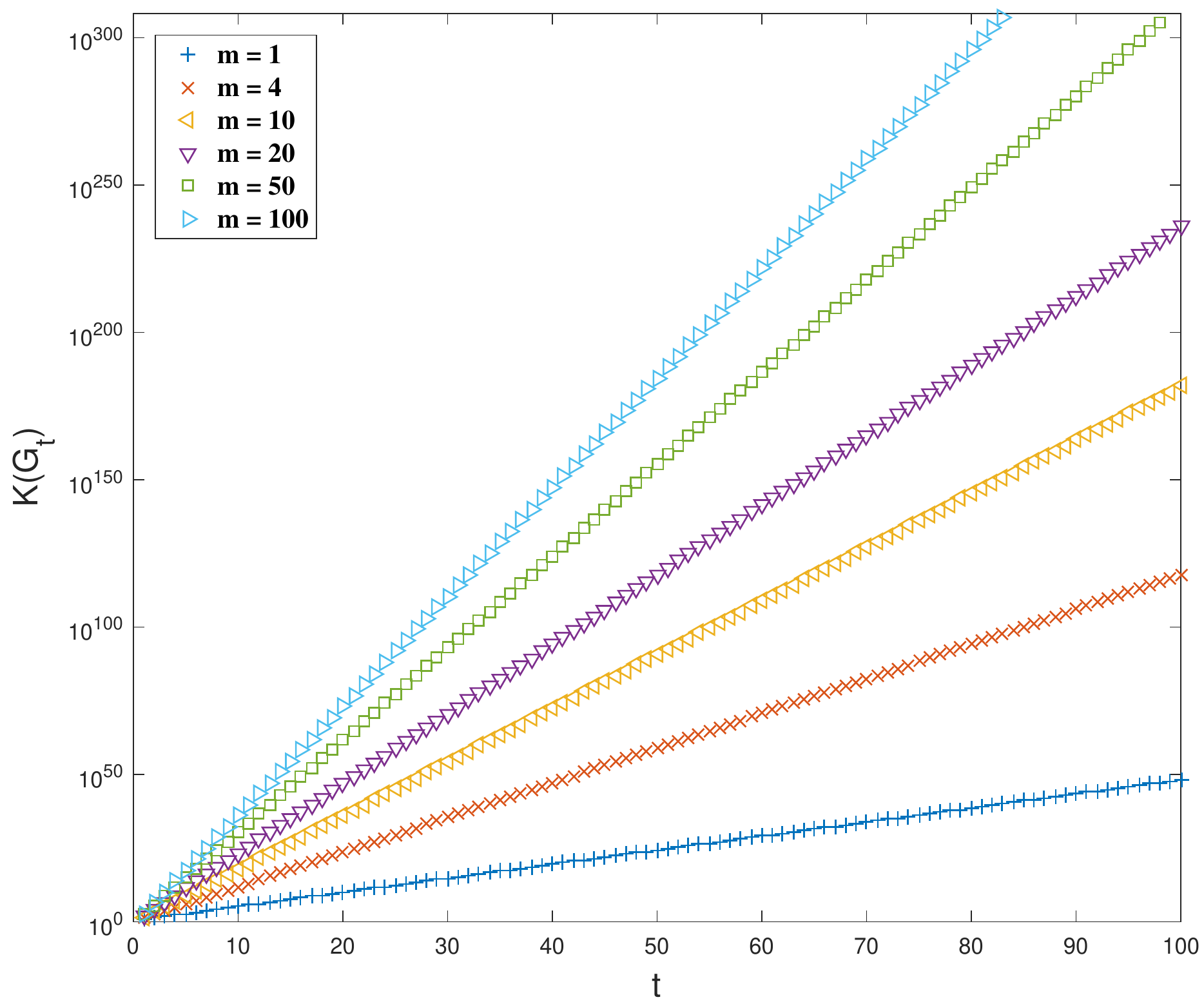}
\caption{Numerical simulation diagram of Kemeny constant $K(G_t)$ with $r=1$.}
\end{minipage}
\begin{minipage}[t]{0.48\textwidth}
\centering
\includegraphics[width=8cm]{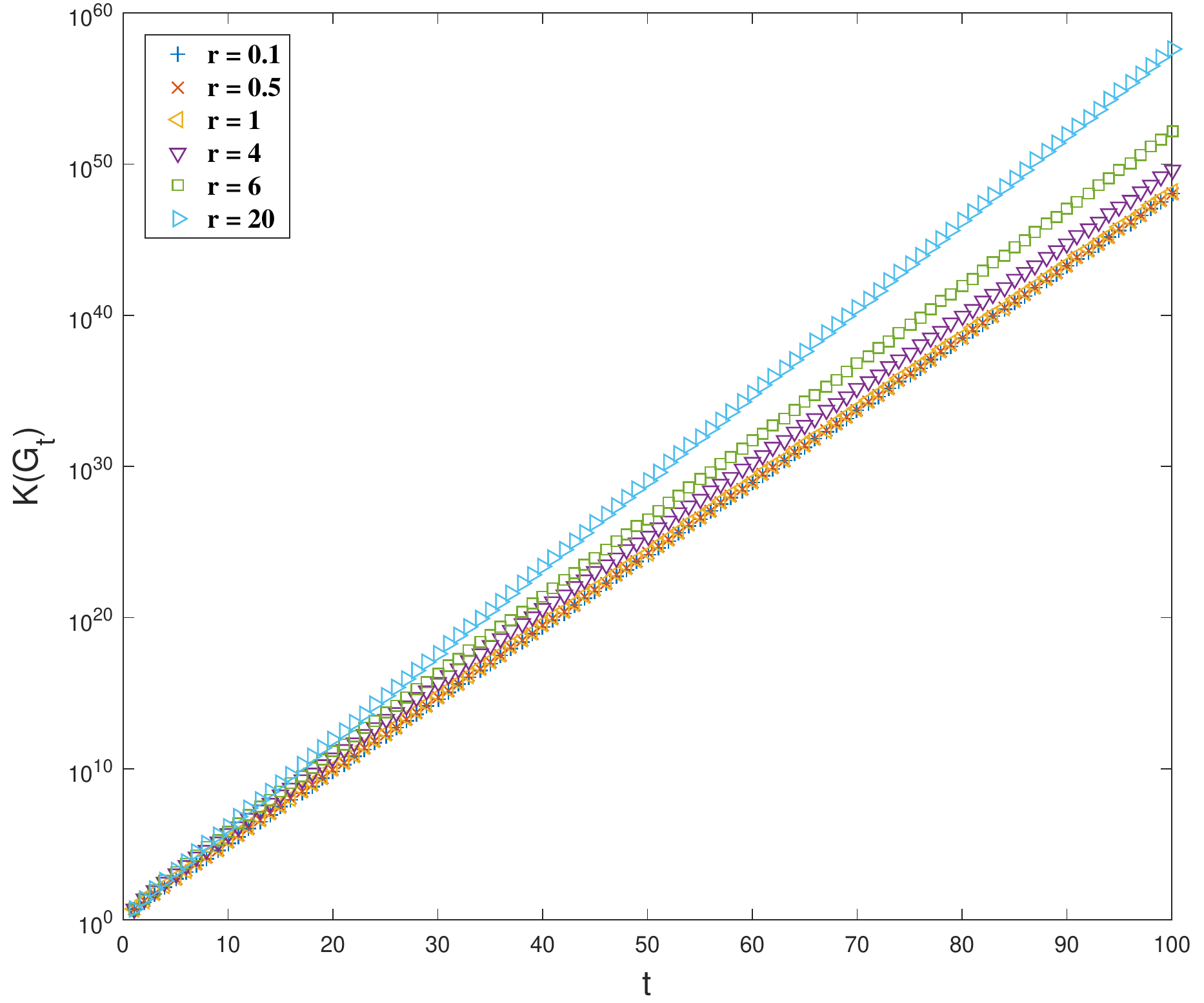}
\caption{Numerical simulation diagram of Kemeny constant $K(G_t)$ with $m=1$.}
\end{minipage}
\end{figure}

\begin{figure}[t]
\centering
\includegraphics[width=12cm]{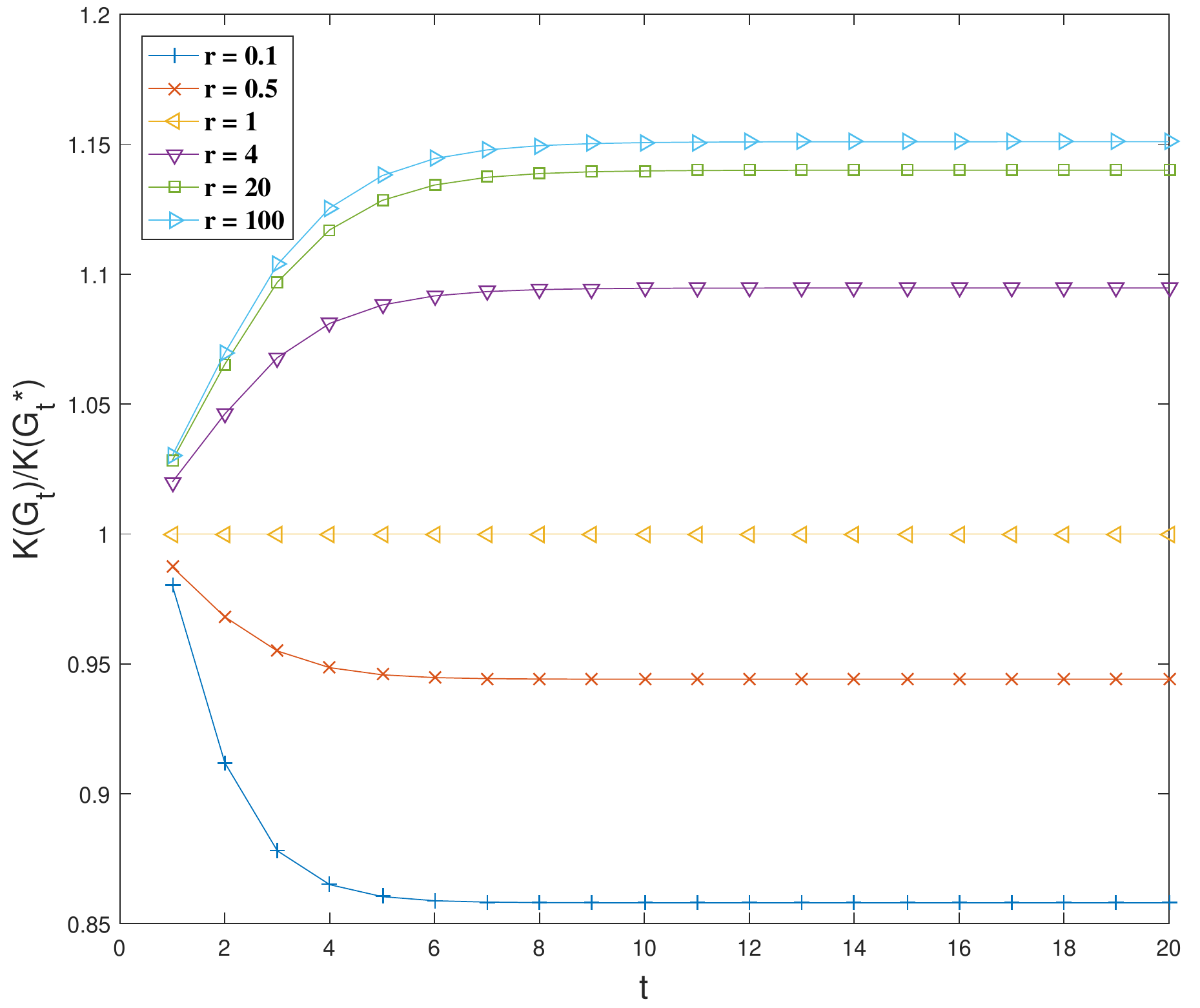}
\caption{Numerical simulation diagram of $K(G_t)/K(G_t^*)$ with $m=4$, where $G_t^*$ represents the network model generated by the network operation $\tau_4^1(\cdot)$.}
\end{figure}

Finally, based on the initial network complete network $K_3$, a numerical simulation of Kemmeny constant is carried out, and the results are shown in Fig.2, Fig.3 and Fig.4. In Fig.2, with the fixed weight factor $r=1$, it can be found that adjusting the scale factor $m$ will have a huge impact on the Kemeny constant of the network. The reason is that the change of scale factor $m$ will have a strong impact on the overall topology of the network. Then, it can be found from Fig.3 that the fixed size factor $m=1$, when $r<4$, Kemeny constant is less affected by the weight factor $r$, but when $r\geq4$, the influence will be significantly increased. According to Eq.(\ref{eq:4.1.1.7}), if and only if $m=1$ and $r\geq4$, the power base in the leading term of Kemeny constant will be affected by the weight factor $r$. Therefore, when $m>1$, the weight factor $r$ affects the Kemmeny constant of the network by changing the coefficient of leading term. In Fig.4, a numerical simulation is carried out for $K(G_t)/K(G_t^*)$ with a fixed size factor $m=4$, where $K(G_t^*)$ represents the Kemeny constant when the weight factor $r=1$. Obviously, by adjusting the weight factor $r$, the overall dynamic properties of the network can be regulated within a certain range without changing the topological structure of the network.

\subsection{The Multiplicative Degree-Kirchhoff Index $\widehat{\mathcal{K}}(G_t)$}

The Kichhoff index is a very important characteristic quantity in the resistance network, and it can be used to measure the overall connectivity of the network.\cite{klein1993resistance,li2018kirchhoff} In recent years, some scholars have improved on the basis of the Kirchhoff index, and proposed the definition of the Multiplicative Degree-Kirchhoff index $\widehat{\mathcal{K}}(G_t)$ on the resistance network $G_t$ as follows:\cite{chen2007resistance}
$$
\widehat{\mathcal{K}}(G_t)=\frac{1}{2}\sum_{i,j=1}^{|V_t|}s_t^is_t^j\cdot r_{ij}=\sum_{\{i,j\}\subseteq V_t}s_t^is_t^j\cdot r_{ij},
$$
where $r_{ij}$ represents the resistance distance between node $i$ and node $j$ in the resistance network $G_t$.
In addition, the Multiplicative Degree-Kirchhoff Index $\widehat{\mathcal{K}}(G_t)$ can also be determined by the spectrum of the normalized adjacency matrix $P_t$, that is, the following relationship is satisfied:
$$
\widehat{\mathcal{K}}(G_t)=2|E_t|\sum_{i=2}^{|V_t|}\frac{1}{1-\lambda_t(i)}=2|E_t|\cdot K(G_t).
$$
Since the analytic expression of Kemeny constant $K(G)$ on the iterative weighted network $G_t$ has been solved, the analytic expression of $\widehat{\mathcal{K}}(G_t)$ naturally satisfies the following relation:

\subsubsection{Case 1: $G_0$ is not a tree}
\begin{eqnarray}
\widehat{\mathcal{K}}(G_t)
=\left\{
\begin{footnotesize}
\begin{array}{ll}
\Big(\widehat{\mathcal{K}}(G_0)-2|E_0|(k_1+k_2)\Big)\big(\frac{[(m+3)mr+2](m^2+3m+2)}{2(mr+2)}\big)^t & \nonumber\\ +2|E_0|k_1\big(\frac{m^2+3m+2}{2}\big)^{2t}+2|E_0|k_2\big(\frac{m^2+3m+2}{2}\big)^t,\nonumber\\ 
\textrm{if $m\neq1$ or $r\neq4$;}\\
&\\
(\widehat{\mathcal{K}}(G_0)+\frac{1}{3}|E_0|^2-\frac{2}{3}|V_0||E_0|+\frac{2}{9}|E_0|)\cdot3^{2t}+\frac{|E_0|^2}{2}3^{2t}\cdot t&\nonumber\\
-(\frac{1}{3}|E_0|^2-\frac{2}{3}|V_0||E_0|+\frac{2}{9}|E_0|)\cdot3^{t}, \nonumber\\
\textrm{if $m=1$ and $r=4$.}
\end{array}
\end{footnotesize}
\right.
\end{eqnarray}
Therefore, as the number of iterations $t$ approaches infinity, the leading term of $\widehat{\mathcal{K}}(G_t)$ satisfy the following approximate relationship:

\begin{eqnarray*}
\widehat{\mathcal{K}}(G_t)
\sim\left\{
\begin{footnotesize}
\begin{array}{ll}
\Big(\widehat{\mathcal{K}}(G_0)-2|E_0|(k_1+k_2)\Big)\big(\frac{[(m+3)mr+2](m^2+3m+2)}{2(mr+2)}\big)^t, \nonumber\\
\textrm{if $m=1$ and $r>4$;}\\
\frac{|E_0|^2}{2}3^{2t}\cdot t, \nonumber\\
\textrm{if $m=1$ and $r=4$};\\
2|E_0|k_1\big(\frac{m^2+3m+2}{2}\big)^{2t},   \nonumber\\
\textrm{others}.
\end{array}
\end{footnotesize}
\right.
\end{eqnarray*}

\subsubsection{Case 2: $G_0$ is a tree}

\begin{eqnarray}
\widehat{\mathcal{K}}(G_t)
=\left\{
\begin{footnotesize}
\begin{array}{ll}
\Big(\widehat{\mathcal{K}}(G_0)+2(|V_0|-1)k_5\Big)\big(\frac{[(m+3)mr+2](m^2+3m+2)}{2(mr+2)}\big)^t&\nonumber\\ +2(|V_0|-1)k_3\big(\frac{m^2+3m+2}{2}\big)^{2t}+2(|V_0|-1)k_4\big(\frac{m^2+3m+2}{2}\big)^t,\nonumber\\
\textrm{if $m\neq1$ or $r\neq4$;}\\
&\\
(\widehat{\mathcal{K}}(G_0)-\frac{29}{72}|V_0|^2+\frac{61}{108}|V_0|-\frac{35}{216})\cdot3^{2t}+\frac{(|V_0|-1)^2}{2}3^{2t}\cdot t&\nonumber\\
+\frac{3|V_0|^2-2|V_0|-1}{8}\cdot3^{t}, \nonumber\\
\textrm{if $m=1$ and $r=4$.}
\end{array}
\end{footnotesize}
\right.
\end{eqnarray}
Naturally, when $t\rightarrow\infty$, the leading term of $\widehat{\mathcal{K}}(G_t)$ obey:
\begin{eqnarray*}
\widehat{\mathcal{K}}(G_t)
\sim\left\{
\begin{footnotesize}
\begin{array}{ll}
\Big(\widehat{\mathcal{K}}(G_0)+2(|V_0|-1)k_5\Big)\big(\frac{[(m+3)mr+2](m^2+3m+2)}{2(mr+2)}\big)^t,\nonumber\\
\textrm{if $m=1$ and $r>4$;}\\
\frac{(|V_0|-1)^2}{2}3^{2t}\cdot t, \nonumber\\
\textrm{if $m=1$ and $r=4$;}\\
2(|V_0|-1)k_3\big(\frac{m^2+3m+2}{2}\big)^{2t}, \nonumber\\
\textrm{others.}\\
\end{array}
\end{footnotesize}
\right.
\end{eqnarray*}

\subsection{The Number of Weighted Spanning Trees}

A spanning tree of a connected network is a minimal connected subgraph containing all nodes in the network.
As a global structure of the network, it is closely related to many aspects of the network, so a study on the spanning tree number of the network can enable us to have a deeper understanding of the influence of the network structure on the process of dynamics. For example, the redundancy of a network can be measured by the number of spanning trees, which is an important indicator to measure the robustness and stability of the network.
Moreover, for a connected network, it is found that the spanning tree has the best synchronization capability when the total number of nodes remains unchanged. In addition, the study of weighted spanning tree can also be used to determine the location of hub nodes in the network.

In weighted iterative network $G_t$, the number of weighted spanning trees is denoted as $N_{tr}^{t}$. Chung et al. have shown that the number of spanning trees in the network can be determined by the non-zero eigenvalues of the normalized Laplacian matrix, which satisfies the following relationship:\cite{chang2014spanning,chung1997spectral}
$$
N_{tr}^t=\frac{\prod_{i=1}^{|V_t|}s_t^i\prod_{i=2}^{|V_t|}\sigma_t^i}{\sum_{i=1}^{|V_t|}s_t^i}.
$$
Obviously, the ratio of the number of weighted spanning trees for two successive generations is:
\begin{eqnarray}\label{eq:4.3.1}
\frac{N_{tr}^t}{N_{tr}^{t-1}}\!&=&\!\frac{\sum_{i=1}^{|V_{t-1}|}\!s_{t-1}^i}{\sum_{i=1}^{|V_t|}\!s_{t}^i}\!\cdot \! \frac{\prod_{i=1}^{|V_t|}\!s_{t}^i}{\prod_{i=1}^{|V_{t-1}|}\!s_{t-1}^i}\!\cdot \! \frac{\prod_{i=2}^{|V_t|}\!\sigma_{t}^i}{\prod_{i=2}^{|V_{t-1}|}\!\sigma_{t-1}^i}.
\end{eqnarray}
The first term on the right side of Eq.(\ref{eq:4.3.1}) obviously satisfies:
\begin{eqnarray}\label{eq:4.3.2}
\frac{\sum_{i=1}^{|V_{t-1}|}\!s_{t-1}^i}{\sum_{i=1}^{|V_t|}\!s_{t}^i}=\frac{S_{t-1}}{S_t}=\frac{2}{m^2r+3mr+2}.
\end{eqnarray}

The second term on the right side of Eq.(\ref{eq:4.3.1}) can be decomposed and expressed as:
\begin{eqnarray}\label{eq:4.3.3}
\frac{\prod_{i=1}^{|V_t|}\!s_{t}^i}{\prod_{i=1}^{|V_{t-1}|}\!s_{t-1}^i}=\prod_{i\in\widehat{V}_t}\frac{s_t^i}{s_{t-1}^i}\cdot
\prod_{i\in\bar{V}_t}s_t^i=(mr+1)^{|V_{t-1}|}\cdot\prod_{i\in\bar{V}_t}s_t^i.
\end{eqnarray}
According to the definition of the weighted $m-$clique annex operation $\tau^r_m(\cdot)$, the weight of edges in the iterated weighted network $G_t$ must be $r^k$ $(0\leq  k\leq t)$, so the number of edges with weight $r^k$ is denoted as $E_t(r^k)$. In addition, it is easy to prove that $E_t(r^k)$ satisfies the following two relationships:
$$
|E_{t}|=\sum_{k=0}^{t}E_{t}(r^k)
$$
and
\begin{eqnarray*}
E_{t}(r^k)
=\left\{
\begin{array}{ll}
E_{t-1}(r^{k})+\frac{m^2+3m}{2}E_{t-1}(r^{k-1}), & 1\leq k \leq t-1 \\
\frac{m^2+3m}{2}E_{t-1}(r^{k-1}), & k=t.
\end{array}
\right.
\end{eqnarray*}
For node $i\in\bar{V}_t$, $s_t^i$ is completely determined by its parent edge $e_i^*$, so the following relationship can be obtained:
\begin{eqnarray}\label{eq:4.3.4}
\prod_{i\in\bar{V}_t}s_t^i&=&\prod_{i\in\bar{V}_t}(m+1)\omega_i^*=\prod_{k=0}^{t-1}\Big[(m+1)r^{k+1}\Big]^{mE_{t-1}(r^{k})}=(m+1)^{m|E_{t-1}|}\cdot r^{m\chi(t-1)},
\end{eqnarray}
where  $\chi(t-1)=\sum_{k=0}^{t-1}(k+1)E_{t-1}(r^k)$.

By decomposing $\chi(t)$, the following iterative relationship can be obtained:
\begin{eqnarray}\label{eq:4.3.5}
\chi(t)&=&\sum_{k=0}^{t}(k+1)E_{t}(r^k)=|E_t|+\sum_{k=1}^{t}kE_{t}(r^k)\nonumber\\
&=&|E_t|+\sum_{k=1}^{t-1}k\Big[E_{t-1}(r^{k})+\frac{m^2+3m}{2}E_{t-1}(r^{k-1})\Big]+ \frac{m^2+3m}{2}tE_{t-1}(r^{t-1})\nonumber\\
&=&|E_t|+\sum_{k=1}^{t-1}kE_{t-1}(r^{k})+\frac{m^2+3m}{2}\sum_{k=0}^{t-1}(k+1)E_{t-1}(r^{k})\nonumber\\
&=&|E_t|-|E_{t-1}|+\frac{m^2+3m+2}{2}\chi(t-1)\nonumber\\
&=&|E_0|\frac{m^2+3m}{2}(\frac{m^2+3m+2}{2})^{t-1}+\frac{m^2+3m+2}{2}\chi(t-1)
\end{eqnarray}
Because of $\chi(0)=|E_0|$, the analytical expression of $\chi(t)$ can be deduced from the iterative relationship in Eq.(\ref{eq:4.3.5}) as follows:
\begin{eqnarray}\label{eq:4.3.6}
\chi(t)\!=\!|E_0|\!\cdot\!\Big[\frac{m^2\!+\!3m\!+\!2}{2}\!+\!\frac{m^2\!+\!3m}{2}t\Big]\!\cdot\!(\frac{m^2\!+\!3m\!+\!2}{2})^{t-1}.
\end{eqnarray}
By substituting Eq.(\ref{eq:4.3.6}) and Eq.(\ref{eq:4.3.4}) into Eq.(\ref{eq:4.3.3}), it can be obtained that:
\begin{eqnarray}\label{eq:4.3.7}
\frac{\prod_{i=1}^{|V_t|}\!s_{t}^i}{\prod_{i=1}^{|V_{t-1}|}\!s_{t-1}^i}&=& (mr+1)^{|V_{t-1}|}\cdot(m+1)^{m|E_{t-1}|}
 \cdot r^{m|E_0|\cdot\big[1+\frac{m^2+3m}{2}t\big](\frac{m^2+3m+2}{2})^{t-2}}.
\end{eqnarray}

Then, when the number of iterations $t>1$ or $G_0$ is not tree, the following relation can be obtained according to the spectral analysis and Eq.(\ref{eq:4.1.4}):
\begin{eqnarray*}
\prod_{i=2}^{|V_t|}\!\sigma_{t}^i&=&(\frac{m+2}{m+1})^{|V_t|-|V_{t-1}|-|E_{t-1}|}\cdot(\frac{2}{m+1})^{|E_{t-1}|-|V_{t-1}|}\cdot \frac{m^2r+3mr+2}{(mr+1)(m+1)}\cdot \prod_{i=2}^{|V_{t-1}|}\big(g_1(\sigma_{t-1}^i)g_2(\sigma_{t-1}^i)\big)\nonumber\\
&=&(\frac{m+2}{m+1})^{|V_t|-|V_{t-1}|-|E_{t-1}|}\cdot(\frac{2}{m+1})^{|E_{t-1}|-|V_{t-1}|}\cdot \frac{m^2r+3mr+2}{(mr+1)(m+1)}\prod_{i=2}^{|V_{t-1}|}\sigma_{t-1}^i\cdot (\frac{mr+2}{(mr+1)(m+1)})^{|V_t|-1}.
\end{eqnarray*}
Therefore, the last term on the right side of Eq.(\ref{eq:4.3.1}) can be written as:
\begin{eqnarray}\label{eq:4.3.8}
\frac{\prod_{i=2}^{|V_t|}\!\sigma_{t}^i}{\prod_{i=2}^{|V_{t-1}|}\!\sigma_{t-1}^i}=(\frac{m\!+\!2}{m\!+\!1})^{|V_t|-|V_{t-1}|-|E_{t-1}|}\!\cdot\! (\frac{2}{m\!+\!1})^{|E_{t-1}|-|V_{t-1}|} \!\cdot\! \frac{m^2r\!+\!3mr\!+\!2}{(mr\!+\!1)(m\!+\!1)}\!\cdot\!(\frac{mr\!+\!2}{(mr\!+\!1)(m\!+\!1)})^{|V_{t-1}|-1}
\end{eqnarray}

By substituting Eq.(\ref{eq:4.3.2}), Eq.(\ref{eq:4.3.7}) and Eq.(\ref{eq:4.3.8}) into Eq.(\ref{eq:4.3.1}), it can be obtained that:
\begin{eqnarray}\label{eq:4.3.9}
\frac{N_{tr}^t}{N_{tr}^{t-1}}&=&r^{m|E_0|\cdot\big[1+\frac{m^2+3m}{2}t\big](\frac{m^2+3m+2}{2})^{t-2}}\cdot (m+2)^{|V_t|-|V_{t-1}|-|E_{t-1}|} \cdot 2^{|E_{t-1}|-|V_{t-1}|+1}\cdot (mr+2)^{|V_{t-1}|-1}
\end{eqnarray}
In addition, when $t=1$ and $G_0$ is tree, according to the above calculation method, it can be obtained that:
\begin{eqnarray*}
\frac{\prod_{i=2}^{|V_1|}\!\sigma_1^i}{\prod_{i=2}^{|V_0|}\!\sigma_0^i}=(\frac{m+2}{m+1})^{|V_1|-2|V_0|}\cdot \frac{m^2r+3mr+2}{(mr+1)(m+1)}\cdot (\frac{mr+2}{(mr+1)(m+1)})^{|V_0|-1}
\end{eqnarray*}
and
\begin{eqnarray}\label{eq:4.3.10}
\frac{N_{tr}^1}{N_{tr}^{0}}=2r^{m(|V_0|-1)}\cdot (m+2)^{|V_1|-2|V_0|}\cdot (mr+2)^{|V_0|-1}
\end{eqnarray}
where it can be noticed that $N_{tr}^{0}=1$.

In view of the influence of different initial network topologies on the number of weighted spanning trees $N_{tr}^t$, we need to conduct a classification discussion on this.

\subsubsection{Case 1: $G_0$ is not a tree}

Under this condition, only the iterative relation of $N_{tr}^t$ in Eq.(\ref{eq:4.3.9}) can be used to obtain following analytic expression:
\begin{eqnarray}
N_{tr}^t&=&\prod_{i=1}^t\frac{N_{tr}^i}{N_{tr}^{i-1}}\cdot N_{tr}^0\nonumber\\
&=& N_{tr}^0\cdot 2^{\frac{|E_0|}{2m(m+3)}[(m+1)(2m+1)^t+4mt-m-1]-(|V_0|-1)t}\cdot r^{\frac{|E_0|}{4m(2m+1)}[(4m^2t+2mt-1)(2m+1)^t+1]}\nonumber\\
&&(mr+2)^{(|V_0|-1)t+\frac{2|E_0|}{m+3}[\frac{(2m+1)^t-1}{2m}-t]}\cdot (m+2)^{\frac{3|E_0|(m-1)}{2m(m+3)}[(2m+1)^t-1]}.
\end{eqnarray}

\subsubsection{Case 2: $G_0$ is a tree}

When the initial network $G_0$ is a tree, Eq.(\ref{eq:4.3.9}) and Eq.(\ref{eq:4.3.10}) need be combined to obtain the analytic expression of $N_{tr}^t$, which can be expressed as follows:
\begin{eqnarray}
N_{tr}^t&=&\prod_{i=2}^t\frac{N_{tr}^i}{N_{tr}^{i-1}}\cdot N_{tr}^1\nonumber\\
&=& 2^{\frac{|V_0|-1}{2m(m+3)}[(m+1)(2m+1)^t+4mt-7m-2m^2-1]-(|V_0|-1)(t-1)+1} \nonumber\\
&&\cdot r^{m(|V_0|-1)[\frac{4m^2t+2mt-1}{4m^2(2m+1)^2}(2m+1)^{t+1}-\frac{4m^2+2m-1}{4m^2}+1]}\nonumber\\
&&\cdot(mr+2)^{(|V_0|-1)[\frac{2}{m+3}(\frac{(2m+1)^t-1}{2m}-t)+t]} \nonumber\\
&&\cdot(m+2)^{\frac{|E_0|-1}{2m(m+3)}[3(m-1)(2m+1)^t+3(m+1)-6m^2]+\frac{4m}{m+3}(|V_0|-1)-|V_0|}.
\end{eqnarray}
Hence, the analytic expressions for the number of the weighted spanning tree $N_{tr}^t$ have been derived in the iterated weighted network $G_t$ generated by the weighted $m-$clique annex operation $\tau^r_m(\cdot)$.

\section{Conclusion}

In this paper, a weighted $m-$clique annex operation $\tau^r_m(\cdot)$ has been proposed and through this operation, a type of weighted iterative network model been constructed. The network model has global small-world property, scale-free property, and clique structures on local features. Then, based on the iterative property of the network structure, the iterative relationship of the eigenvalues of the normalized Laplacian matrix is analyzed, and three spectral applications are given, namely, Kenemy constant, Multiplicative Degree-Kirchhoff Index and the number of weighted spanning trees.

In addition, because Kemeny constant is closely related to the first-passage time of the network, which can reflect the functions and dynamic properties of the network remarkably, this characteristic quantity is taken as an example to analyze the impact of size factor $m$ and weight factor $r$ on network model. The results show that under the same number of iterations, by changing the size factor $m$, the topology of the network can be significantly changed, and the overall function and dynamic properties of the network can be dramatically affected, which is consistent with the inference in Remark 2. In comparison, adjusting the weight factor $r$ can only control the dynamic properties and functions of the network in a small range, but its advantage is that the structure of the network has not changed, so the cost of regulation is lower. Spontaneously, the simultaneous regulation of the above two factors can make the function and dynamic properties of the network model reach the ideal state. Given the wide application of the first-passage time and Kemeny constant in measuring navigation efficiency, random search cost and  efficiency of robotic surveillance in the network,\cite{levene2002kemeny,patel2015robotic} we believe that this adjustable weighted $m-$clique annex operation $\tau^r_m(\cdot)$ and the complex network model generated by it have great application potential in the field of artificial networks.

\section*{Acknowledgments}
This work was supported by Key Project of Natural Science Foundation of China (No. 61833005) and National Natural Science Foundation of China (No. 12026214).

\section*{Data Availability Statement}
The data that support the findings of this study are available within the article.

\bibliographystyle{ws-ijbc}
\bibliography{Ref}
\end{document}